\documentclass[a4paper,11pt]{amsart}

\usepackage{graphicx}
\usepackage{mathptmx}
\usepackage{amsmath}
\usepackage{amssymb}
\usepackage{enumitem}
\usepackage{xcolor}

\newmuskip\pFqmuskip

\newcommand*\pFq[6][8]{%
  \begingroup % only local assignments
  \pFqmuskip=#1mu\relax
  \mathcode`=\string"8000
  \begingroup\lccode`\~=`\,
  \lowercase{\endgroup\let~}\pFqcomma
  F^{#2}_{#3}{\left(\genfrac..{0pt}{}{#4}{#5}\bigg|#6\right)}%
  \endgroup
}
\newcommand{\pFqcomma}{\mskip\pFqmuskip}

\newtheorem{theorem}{Theorem}[section]

\newtheorem{corollary}[theorem]{Corollary}

\begin{document}

\title[]{Explicit formulas related to Euler's product expansion for cosine function}

\author{Taekyun  Kim}
\address{Department of Mathematics, Kwangwoon University, Seoul 139-701, Republic of Korea}
\email{tkkim@kw.ac.kr}
\author{Dae San  Kim}
\address{Department of Mathematics, Sogang University, Seoul 121-742, Republic of Korea}
\email{dskim@sogang.ac.kr}


\begin{abstract}
In this paper, we derive by using elementary methods some continued fractions, certain identities involving derivatives of $\tan x$, several expressions for $\log \cosh x$ and an identity for $\pi^{2}$, from a series expansion of $\tan x$, which gives the product expansion of the cosine function.
\end{abstract}

\maketitle

\markboth{\centerline{\scriptsize Explicit formulas related to Euler's product expansion for cosine function}}
{\centerline{\scriptsize T. Kim and D. S. Kim}}

\section{Introduction and preliminaries}
In this paper, from a series expansion of $\tan x$ (see \eqref{0}), which gives the product expansion of the cosine function, we derive some continued fractions, certain identities involving derivatives of $\tan x$, several expressions for $\log \cosh x$ and an identity for $\pi^{2}$. \par
In more detail, the outline of this paper is as follows. In Section 1, we recall the Bernoulli numbers, the Riemann zeta function together with their special values at even positive integers, and the product expansion of the sine function. We remind the reader of the Euler numbers and explicit expressions of them obtained by a contour integral. In fact, these yield a relationship between Euler and Bernoulli numbers (see \eqref{11}). We show that
\begin{align}
\tan x&=\sum_{k=0}^{\infty}\bigg(\frac{2}{\pi(2k+1)}\bigg)^{2}2x\frac{1}{1-\big(\frac{2x}{\pi(2k+1)}\big)^{2}} \label{0} \\
&=8x\bigg(\frac{1}{\pi^{2}-4x^{2}}+\sum_{k=1}^{\infty}\frac{1}{((2k+1)\pi)^{2}-4x^{2}}\bigg),\nonumber
\end{align}
by expressing $\tan x$ as a series involving the Euler numbers and invoking the aforementioned explicit expressions of Euler numbers. Then we derive the product expansion of the cosine function from \eqref{0}. Section 2 contains the main results of this paper. In Theorem 2.1, we derive a continued fraction for $\frac{\pi}{4}$. In Theorem 2.2, by successively applying differential operators to the expression $\frac{\tan x}{8x}+\frac{1}{4x^2-\pi^2}=\sum_{k=1}^{\infty}\frac{1}{((2k+1)\pi)^2-4x^2}$ we show that it is equal to $\sum_{k=1}^{\infty}\Big(\frac{1}{((2k+1)\pi)^2-4x^2}\Big)^{n+1}$. By taking $(2n-2)$th derivative of the partial fraction expression, $\tan x=2\sum_{k=0}^{\infty}\Big(\frac{1}{(2k+1)\pi -2x}-\frac{1}{(2k-1)\pi +2x}\Big)$, of \eqref{0}, dividing by $2^{2n-1}(2n-2)!$, and evaluating at $x=\frac{\pi}{4}$, we have the expression in Theorem 2.3, which is equal to $\frac{1}{4^{2n-1}}\Big(\zeta\Big(2n-1,\frac{1}{4}\Big)-\zeta\Big(2n-1,\frac{3}{4}\Big)\Big)$. Here $\zeta(s,x)$ is the Hurwitz zeta function. By taking the derivative of \eqref{0}, and letting $x=\frac{\pi}{2n+1}$ we obtain an expression of $\pi^{2}$ as an infinite series, which is valid for any integer $n$. Letting $n=0$, we get an expression of $\frac{\pi}{2}$ in Corollary 2.5. From the product expansion of the cosine function, we show that $\frac{\tanh x-\tan x}{8x}=\sum_{k=0}^{\infty}\bigg(\frac{1}{(2k+1)^{2}\pi^{2}+4x^{2}}-\frac{1}{(2k+1)^{2}\pi^{2}-4x^{2}}\bigg)$. From the identity obtained by letting $x=\frac{\pi}{2(2n+1)}$ in this, we get a continued fraction for $\bigg(\frac{\pi}{4}\frac{\tanh(\frac{\pi}{2(2n+1)})-\tan(\frac{\pi}{2(2n+1)})}{2n+1}\bigg)^{-1}$ in Theorem 2.6. Finally, from the observation $\tanh x=-\sum_{n=1}^{\infty}E_{2n-1}\frac{2^{2n-1}x^{2n-1}}{(2n-1)!}$ we have several expressions for $\log \cosh x$ in Theorem 2.7. For the rest of this section, we recall the facts that are needed throughout this paper.

\vspace{0.1in}

The Bernoulli numbers are defined by
\begin{equation}
\frac{z}{e^{z}-1}=\sum_{n=0}^{\infty}B_{n}\frac{z^{n}}{n!},\quad (\mathrm{see}\ [1-20]).\label{1}
\end{equation}
The first few terms of $B_n$ are given by:
\begin{align*}
&B_0=1,\,B_1=-\frac{1}{2},\,B_2=\frac{1}{6},\,B_4=-\frac{1}{30},\,B_6=\frac{1}{42},\,B_8=-\frac{1}{30},\,B_{10}=\frac{5}{66},\,\\
&B_{12}=-\frac{691}{2730},\, B_{14}=\frac{7}{6},\,B_{16}=-\frac{3617}{510},\,B_{18}=\frac{43867}{798},\,B_{20}=-\frac{174611}{330},\dots; \\
&B_{2k+1}=0,\,\,(k \ge 1).
\end{align*} \par
For $s\in\mathbb{C}$ with $\mathrm{Re}(s)>1$, the Riemann zeta function is defined by
\begin{equation}
\zeta(s)=\sum_{n=1}^{\infty}\frac{1}{n^{s}},\quad (\mathrm{see}\ [1,9,12]). \label{2}	
\end{equation}
For $n\in\mathbb{N}$, we have
\begin{equation}
\zeta(2n)=\sum_{k=1}^{\infty}\frac{1}{k^{2n}}=\frac{(-1)^{n-1}(2\pi)^{2n}}{2(2n)!}B_{2n},\quad (\mathrm{see}\ [10]).\label{3}
\end{equation}
We observe from \eqref{1} that
\begin{equation}
z\cot z=\frac{2iz}{e^{2iz}-1}+iz=1+\sum_{n=1}^{\infty}\frac{(-1)^{n}2^{2n}B_{2n}}{(2n)!}z^{2n}.\label{4}
\end{equation}
From \eqref{3} and \eqref{4}, we get the following product expansion of the sine function:
\begin{equation*}
\frac{\sin z}{z}=\prod_{k=1}^{\infty}\bigg(1-\bigg(\frac{z}{k\pi}\bigg)^{2}\bigg),\quad (\mathrm{see}\ [2,10,13]).
\end{equation*} \par
It is well known that a continued fraction for $\pi$ is given by
\begin{equation}
\pi=\cfrac{4}{1+\cfrac{1^{2}}{2+\cfrac{3^{2}}{2+\cfrac{5^{2}}{2+\cfrac{7^{2}}{2+\cfrac{9^{2}}{2+\cdots } } } } } },\quad (\mathrm{see}\ [5,14]), \label{5}
\end{equation}
In [14], another continued fraction for $\pi$ is given by
\begin{equation}
	\pi=3+\cfrac{1^{2}}{6+\cfrac{3^{2}}{6+\cfrac{5^{2}}{6+\cfrac{7^{2}}{6+\cfrac{9^{2}}{6+\cdots} } } } }\label{6}
\end{equation} \par
The Euler numbers are defined by
\begin{equation}
\frac{2}{e^{z}+1}=\sum_{n=0}^{\infty}E_{n}\frac{z^{n}}{n!},\quad (\mathrm{see}\ [8,10,11]).\label{7}
\end{equation}
From \eqref{7}, we note that
\begin{equation}
E_{0}=1,\quad E_{n}=-\sum_{i=0}^{n}\binom{n}{i}E_{i},\quad (n\in\mathbb{N}),\quad (\mathrm{see}\ [10]).\label{8}
\end{equation}
The first few terms of $E_n$ are given by:
\begin{align*}
&E_0=1,\,E_1=-\frac{1}{2},\,E_3=\frac{1}{4},\,E_5=-\frac{1}{2},\,E_7=\frac{17}{8},\,E_9=-\frac{31}{2},\,E_{11}=\frac{691}{4},\,\\
&E_{13}=-\frac{5461}{2},\, E_{15}=\frac{929569}{16},\,E_{17}=-\frac{3202291}{2},\,E_{19}=\frac{221930581}{4},\dots;\\
& E_{2k}=0,\,\,(k \ge 1).
\end{align*}

For $n\ge 2$, by \eqref{7}, we get
\begin{align*}
0&=\lim_{N \rightarrow \infty}\oint_{C_N}\frac{2}{e^{z}+1}\frac{1}{z^{n}} dz=2 \pi i \Big(\mathop{\mathrm{Res}}_{z=0}\frac{2}{e^{z}+1}\frac{1}{z^{n}}+\sum_{k=-\infty}^{\infty}	\mathop{\mathrm{Res}}_{z=(2k+1)i\pi}\frac{2}{e^{z}+1}\frac{1}{z^{n}}\Big) \\
&=2 \pi i \Big(\frac{E_{n-1}}{(n-1)!}+\sum_{k=-\infty}^{\infty}\frac{-2}{((2k+1)i \pi)^{n}} \Big),
\end{align*}
where the contour $C_{N}$ is the (positively oriented) circle with radius $2N \pi$ centered at the origin $(N=1,2,3, \dots)$. Thus we have
\begin{align}
E_{n-1}&=2(n-1)!\frac{1}{(\pi i)^{n}}\sum_{k=-\infty}^{\infty}\frac{1}{(2k+1)^{n}} \label{9} \\
&=\left\{\begin{array}{ccc}
	\displaystyle 4(n-1)!\frac{(-1)^{n/2}}{\pi^{n}}\sum_{k=0}^{\infty}\frac{1}{(2k+1)^{n}}\displaystyle, & \textrm{for $n$ even,}\\
	0, & \textrm{for $n$ odd.}
\end{array}\right.\nonumber
\end{align}
From \eqref{9}, we note that
\begin{align}
E_{2n-1}&=4(2n-1)!\frac{(-1)^{n}}{\pi^{2n}}\sum_{k=0}^{\infty}\frac{1}{(2k+1)^{2n}},
\label{10} \\
&=4(2n-1)!\frac{(-1)^{n}}{\pi^{2n}}\Big(1-\frac{1}{4^{n}}\Big)\zeta(2n) \nonumber \\
&=-\frac{2^{2n}-1}{n}B_{2n},\quad (n\ge 1). \nonumber
\end{align}
As $E_{2n}=B_{2n+1}=0$\, for $n \ge 1$, the following holds true:
\begin{equation}
E_{n}=-\frac{2(2^{n+1}-1)}{n+1}B_{n+1}, \quad (n \ge 0). \label{11}
\end{equation} \par
We observe that
\begin{align}
i\tan x&=i\frac{\frac{e^{ix}-e^{-ix}}{2i}}{\frac{e^{ix}+e^{-ix}}{2}}=1-\frac{2}{e^{2ix}+1}\label{12}\\
&=1-\sum_{n=0}^{\infty}E_{n}\frac{(2i)^{n}}{n!}x^{n}\nonumber \\
&=-i\sum_{n=0}^{\infty}E_{2n+1}\frac{2^{2n+1}(-1)^{n}}{(2n+1)!}x^{2n+1}.\nonumber	
\end{align}
From \eqref{10} and \eqref{12}, we have
\begin{align}
\tan x&=\sum_{n=0}^{\infty}\frac{2^{2n+1}x^{2n+1}}{(2n+1)!}(-1)^{n-1}E_{2n+1} \label{13} \\
&=\sum_{n=0}^{\infty}\frac{2^{2n+1}x^{2n+1}}{(2n+1)!}\frac{4(2n+1)!}{\pi^{2n+2}}\sum_{k=0}^{\infty}\frac{1}{(2k+1)^{2n+2}}\nonumber\\
&=\sum_{k=0}^{\infty}\bigg(\frac{2}{\pi(2k+1)}\bigg)^{2}2x\sum_{n=0}^{\infty}\bigg(\frac{2x}{\pi(2k+1)}\bigg)^{2n}\nonumber\\
&=\sum_{k=0}^{\infty}\bigg(\frac{2}{\pi(2k+1)}\bigg)^{2}2x\frac{1}{1-\big(\frac{2x}{\pi(2k+1)}\big)^{2}}.\nonumber	
\end{align}
By \eqref{13}, we get
\begin{align}
-\frac{d}{dx}\log\cos x=\tan x&=-\sum_{k=0}^{\infty}\bigg(\frac{2}{(2k+1)\pi}\bigg)^{2}(-2x)\frac{1}{1-\big(\frac{2x}{(2k+1)\pi}\big)^{2}}\label{14} \\
&=-\sum_{k=0}^{\infty}\frac{d}{dx}\log\bigg(1-\bigg(\frac{2x}{(2k+1)\pi}\bigg)^{2}\bigg).\nonumber
\end{align}
From \eqref{14}, we have
\begin{equation}
\begin{aligned}
	\log\cos x&=\sum_{k=0}^{\infty}\log\bigg(1-\bigg(\frac{x}{(2k+1)\pi}\bigg)^{2}\bigg)\\
	&=\log\prod_{k=0}^{\infty}\bigg(1-\bigg(\frac{2x}{(2k+1)\pi}\bigg)^{2}\bigg).
\end{aligned}\label{15}	
\end{equation}
Hence, by \eqref{15}, we get the product expansion of the cosine function:
\begin{equation}
\cos x=\prod_{k=0}^{\infty}\bigg(1-\bigg(\frac{2x}{(2k+1)\pi}\bigg)^{2}\bigg).\label{16}	
\end{equation}

\section{Explicit formulas related to Euler product expansion for cosine function}
From \eqref{13}, we note that
\begin{align}
\tan x&=\sum_{k=0}^{\infty}\bigg(\frac{2}{(2k+1)\pi}\bigg)^{2}2x\bigg(\frac{1}{1-\big(\frac{2x}{(2k+1)\pi}\big)^{2}}\bigg)\label{17} \\
&=8x\bigg(\frac{1}{\pi^{2}-4x^{2}}+\sum_{k=1}^{\infty}\frac{1}{((2k+1)\pi)^{2}-4x^{2}}\bigg).\nonumber
\end{align}
From \eqref{17}, we have
\begin{equation}
\frac{1}{8x}\tan x=\frac{1}{\pi^{2}-4x^{2}}+\sum_{k=1}^{\infty}\frac{1}{((2k+1)\pi)^{2}-4x^{2}}.\label{18}
\end{equation}
Let $x=\pi$ in \eqref{18}. Then we have
\begin{equation}
\begin{aligned}
\frac{1}{3\pi^{2}}&=\sum_{k=1}^{\infty}\frac{1}{((2k+1)\pi)^{2}-4\pi^{2}}\\
&=\sum_{k=1}^{\infty}\frac{1}{(2k-1)(2k+3)\pi^{2}}.
\end{aligned} \label{19}	
\end{equation}
Multiplying $\pi^{2}$ on both sides of \eqref{19}, we get
\begin{equation}
\frac{1}{3}=\sum_{k=1}^{\infty}\frac{1}{(2k-1)(2k+3)}=\frac{1}{1\cdot 5}+\frac{1}{3\cdot 7}+\frac{1}{5\cdot 9}+\cdots	\label{20}
\end{equation}
Let $x=\frac{\pi}{4}$ in \eqref{18}. Then we see that
\begin{equation}
\frac{1}{2\pi}-\frac{4}{3\pi^{2}}=\sum_{k=1}^{\infty}\frac{4}{(4k+1)(4k+3)\pi^{2}}. \label{21}
\end{equation}
Multiplying $\frac{\pi^{2}}{4}$ on both sides of \eqref{21}, we get
\begin{equation}
\frac{\pi}{8}-\frac{1}{3}=\frac{1}{2}\sum_{k=1}^{\infty}\bigg(\frac{1}{4k+1}-\frac{1}{4k+3}\bigg). \label{22}
\end{equation}
Thus, by \eqref{22}, we get
\begin{equation}
\frac{3\pi-8}{12}=\sum_{k=1}^{\infty}\bigg(\frac{1}{4k+1}-\frac{1}{4k+3}\bigg)=\frac{1}{5}-\frac{1}{7}+\frac{1}{9}-\frac{1}{11}+\frac{1}{13}-\frac{1}{15}+\cdots. \label{23}
\end{equation} \par
From \eqref{20}, we get
\begin{align}
&\frac{12}{3\pi-8}-1=\cfrac{1}{\cfrac{3\pi-8}{12}}-1=\cfrac{1-\frac{1}{5}+\frac{1}{7}-\frac{1}{9}+\frac{1}{11}-\frac{1}{13}+\frac{1}{15}-\cdots}{\frac{1}{5}-\frac{1}{7}+\frac{1}{9}-\frac{1}{11}+\frac{1}{13}-\frac{1}{15}+\cdots}\label{24} \\
&=\cfrac{\big(\frac{4}{5}-\frac{4}{7}+\frac{4}{9}-\frac{4}{11}+\frac{4}{13}-\frac{4}{15}+\cdots)+\big(\frac{5}{7}-\frac{5}{9}+\frac{5}{11}-\frac{5}{13}+\cdots\big)}{\frac{1}{5}-\frac{1}{7}+\frac{1}{9}-\frac{1}{11}+\frac{1}{13}-\frac{1}{15}+\cdots}\nonumber\\
&=4+\cfrac{5^{2}}{\cfrac{\big(1-\frac{5}{7}\big)+\frac{5}{9}-\frac{5}{11}+\frac{5}{13}-\frac{5}{15}+\cdots}{\frac{1}{7}-\frac{1}{9}+\frac{1}{11}-\frac{1}{13}+\frac{1}{15}-\cdots}}\nonumber \\
&=4+\cfrac{5^{2}}{\cfrac{\big(\frac{2}{7}-\frac{2}{9}+\frac{2}{11}-\frac{2}{13}+\cdots\big)+\big(\frac{7}{9}-\frac{7}{11}+\frac{7}{13}-\cdots\big)}{\frac{1}{7}-\frac{1}{9}+\frac{1}{11}-\frac{1}{13}+\frac{1}{15}-\cdots}}\nonumber\\
&=4+\cfrac{5^{2}}{2+\cfrac{7^{2}}{\cfrac{\big(1-\frac{7}{9}\big)+\frac{7}{11}-\frac{7}{13}+\frac{7}{15}-\cdots}{\frac{1}{9}-\frac{1}{11}+\frac{1}{13}-\frac{1}{15}+\cdots}}}\nonumber\\
&=\cdots\nonumber \\
&=4+ \cfrac{5^{2}}{2+\cfrac{7^{2}}{2+\cfrac{9^{2}}{2+\cfrac{11^{2}}{2+\cfrac{13^{2}}{2+\cdots} } } } }.\nonumber
\end{align}
Therefore, by \eqref{24}, we obtain the following theorem.
\begin{theorem}
We have the following continued fraction for $\frac{\pi}{4}$:
\begin{displaymath}
\frac{\pi}{4}=\frac{2}{3}+\cfrac{1}{5+\cfrac{5^{2}}{2+\cfrac{7^{2}}{2+\cfrac{9^{2}}{2+\cfrac{11^{2}}{2+\cdots} } } } }
\end{displaymath}
\end{theorem}
By \eqref{18}, we get
\begin{align*}
	&\frac{1}{8x}\frac{d}{dx}\bigg(\frac{\tan x}{8x}+\frac{1}{4x^{2}-\pi^{2}}\bigg)=\sum_{k=1}^{\infty}\frac{1}{\big(((2k+1)\pi)^{2}-4x^{2}\big)^{2}}, \\
	&\frac{1}{16x}\frac{d}{dx}\frac{1}{8x}\frac{d}{dx}\bigg(\frac{\tan x}{8x}+\frac{1}{4x^{2}-\pi^{2}}\bigg)=\sum_{k=1}^{\infty}\frac{1}{\big(((2k+1)\pi)^{2}-4x^{2}\big)^{3}}.
\end{align*}
Continuing this process, we have
\begin{equation}
\begin{aligned}
	&\sum_{k=1}^{\infty}\bigg(\frac{1}{((2k+1)\pi)^{2}-4x^{2}}\bigg)^{n+1}\\
	&=\frac{1}{8nx}\frac{d}{dx}\frac{1}{8(n-1)x}\frac{d}{dx}\cdots\frac{1}{16x}\frac{d}{dx}\frac{1}{8x}\frac{d}{dx}\bigg(\frac{\tan x}{8x}+\frac{1}{4x^{2}-\pi^{2}}\bigg), \label{26}
\end{aligned}
\end{equation}
where $n$ is a positive integer. \\
From \eqref{26}, we note that
\begin{equation*}
\begin{aligned}
	&\bigg[\frac{1}{8nx}\frac{d}{dx}\frac{1}{8(n-1)x}\frac{d}{dx}\cdots\frac{1}{16x}\frac{d}{dx}\frac{1}{8x}\frac{d}{dx}\bigg(\frac{\tan x}{8x}+\frac{1}{4x^{2}-\pi^{2}}\bigg)\bigg]_{x=\frac{\pi}{4}}\\
	&\quad =\frac{4^{n+1}}{\pi^{2n+2}}\sum_{k=1}^{\infty}\bigg(\frac{1}{(4k+1)(4k+3)}\bigg)^{n+1}.
	\end{aligned}
\end{equation*}
\begin{theorem}
For $n\in\mathbb{N}$, we have
\begin{align*}
&\sum_{k=1}^{\infty}\bigg(\frac{1}{((2k+1)\pi)^{2}-4x^{2}}\bigg)^{n+1}\\
&=\frac{1}{8nx}\frac{d}{dx}\frac{1}{8(n-1)x}\frac{d}{dx}\cdots\frac{1}{16x}\frac{d}{dx}\frac{1}{8x}\frac{d}{dx}\bigg(\frac{\tan x}{8x}+\frac{1}{4x^{2}-\pi^{2}}\bigg).
\end{align*}
\end{theorem}
By \eqref{17}, we get
\begin{align}
\tan x &=\sum_{k=0}^{\infty}\frac{8x}{\big((2k+1)\pi\big)^{2}-4x^{2}}=\sum_{k=0}^{\infty}\frac{8x}{\big((2k+1)\pi-2x)((2k-1)\pi +2x)}\label{27}\\
&=2\sum_{k=0}^{\infty}\bigg(\frac{1}{(2k+1)\pi -2x}-\frac{1}{(2k-1)\pi +2x}\bigg).\nonumber
\end{align}
From \eqref{27}, we have
\begin{align}
&\frac{1}{2^{2n-1}}\frac{1}{(2n-2)!}\frac{d^{2n-2}}{dx^{2n-2}}\tan x \label{28} \\
&=\sum_{k=0}^{\infty}\bigg(\frac{1}{\big((2k+1)\pi-2x\big)^{2n-1}}-\frac{1}{\big((2k+1)\pi+2x\big)^{2n-1}}\bigg), \nonumber
\end{align}
where $n\in\mathbb{N}$.
Taking $n=1$ and $x=\frac{\pi}{4}$ in \eqref{28}, we get
\begin{equation*}
	\frac{\pi}{4}=\sum_{k=0}^{\infty}\bigg(\frac{1}{4k+1}-\frac{1}{4k+3}\bigg)=1-\frac{1}{3}+\frac{1}{5}-\frac{1}{7}+\frac{1}{9}-\frac{1}{11}+\cdots,
\end{equation*}
which agrees with the result in \eqref{23}. \\
Letting $x=\frac{\pi}{4}$ in \eqref{28}, we have
\begin{equation}
\frac{1}{2^{2n-1}}\frac{1}{(2n-2)!}\bigg[\frac{d^{2n-2}}{dx^{2n-2}}\tan x\bigg]_{x=\frac{\pi}{4}}=\sum_{k=0}^{\infty}\bigg(\frac{1}{(4k+1)^{2n-1}}-\frac{1}{(4k+3)^{2n-1}}\bigg)\bigg(\frac{2}{\pi}\bigg)^{2n-1}.\label{31}	
\end{equation}
Thus, by \eqref{31}, we obtain the following theorem.
\begin{theorem}
For $n\ge 2$, we have
\begin{align*}
\bigg(\frac{\pi}{4}\bigg)^{2n-1}\frac{1}{(2n-2)!}\bigg[\frac{d^{2n-2}}{dx^{2n-2}}\tan x\bigg]_{x=\frac{\pi}{4}}&=\sum_{k=0}^{\infty}\frac{1}{(4k+1)^{2n-1}}-\sum_{k=0}^{\infty}\frac{1}{(4k+3)^{2n-1}} \\
&=\frac{1}{4^{2n-1}}\Big(\zeta\Big(2n-1,\frac{1}{4}\Big)-\zeta\Big(2n-1,\frac{3}{4}\Big)\Big),
\end{align*}
where $\zeta(s,x)$ is the Hurwitz zeta function given by
\begin{displaymath}
\zeta(s,x)=\sum_{n=0}^{\infty}\frac{1}{(n+x)^{s}},\quad s\in\mathbb{C}\ \ \mathrm{with}\quad \mathrm{Re}(s)>1.
\end{displaymath}
\end{theorem}
By \eqref{17}, we get
\begin{align}
\sec^{2}x&=\frac{d}{dx}\tan x=8\sum_{k=0}^{\infty}\frac{1}{\big((2k+1)\pi\big)^{2}-(2x)^{2}}+8x\sum_{k=0}^{\infty}\frac{8x}{\big((2k+1)\pi\big)^{2}-(2x)^{2}\big)^{2}}\label{32}\\
&=8\sum_{k=0}^{\infty}\frac{\big((2k+1)\pi\big)^{2}+4x^{2}}{\big(\big((2k+1)\pi\big)^{2}-4x^{2}\big)^{2}}.\nonumber
\end{align}
From \eqref{32}, we note that
\begin{equation}
\frac{1}{8}\sec^{2}x=\sum_{k=0}^{\infty}\frac{(2k+1)^{2}\pi^{2}+4x^{2}}{\big(4x^{2}-(2k+1)^{2}\pi^{2}\big)^{2}}. \label{33}
\end{equation}
Let $x=\frac{\pi}{2n+1},\ (n\in\mathbb{Z})$ in \eqref{33}. Then we have
\begin{align}
\frac{1}{8}\sec^{2}\bigg(\frac{\pi}{2n+1}\bigg)&=\sum_{k=0}^{\infty}\frac{(2k+1)^{2}\pi^{2}+4\big(\frac{\pi}{2n+1}\big)^{2}}{\big(4\big(\frac{\pi}{2n+1}\big)^{2}-(2k+1)^{2}\pi^{2}\big)^{2}} \label{34} \\
&=\sum_{k=0}^{\infty}\frac{\big(4+(2n+1)^{2}(2k+1)^{2}\big) \big(\frac{\pi}{2n+1}\big)^{2}}{\big(4-(2n+1)^{2}(2k+1)^{2}\big)^{2} \big(\frac{\pi}{2n+1}\big)^{4}}.\nonumber
\end{align}
Thus, by \eqref{34}, we get the following theorem.
\begin{theorem}
For $n\in\mathbb{Z}$, we have
\begin{displaymath}
\pi^{2}=8(2n+1)^{2}\cos^{2}\bigg(\frac{\pi}{2n+1}\bigg)\sum_{k=0}^{\infty}\frac{4+(2n+1)^{2}(2k+1)^{2}}{\big(4-(2n+1)^{2}(2k+1)^{2}\big)^{2}}.
\end{displaymath}
\end{theorem}
Let $n=0$ in Theorem 2.5. Then we have
\begin{equation}
\pi^{2}=8\sum_{k=0}^{\infty}\frac{4+(2k+1)^{2}}{\big(4-(2k+1)^{2}\big)^{2}}. \label{36}
\end{equation}
Thus, by \eqref{36}, we obtain the following corollary.
\begin{corollary}
\begin{displaymath}
\frac{\pi}{2} =\sqrt{2\sum_{k=0}^{\infty}\frac{4+(2k+1)^{2}}{\big(4-(2k+1)^{2}\big)^{2}}}.
\end{displaymath}
\end{corollary}
From \eqref{16}, we note that
\begin{equation}
\cos ix=\prod_{k=0}^{\infty}\bigg(1+\frac{4x^{2}}{(2k+1)^{2}\pi^{2}}\bigg).
\label{37}
\end{equation}
As $\cos ix=\cosh x$, from \eqref{37} we have
\begin{equation}
\cosh x=\prod_{k=0}^{\infty}\bigg(1+\bigg(\frac{2x}{(2k+1)\pi}\bigg)^{2}\bigg). \label{39}
\end{equation}
By \eqref{16} and \eqref{39}, we get
\begin{align}
\cosh x\cos x&=\prod_{k=0}^{\infty}\bigg(1+\bigg(\frac{2x}{(2k+1)\pi}\bigg)^{2}\bigg)\prod_{k=0}^{\infty}\bigg(1-\bigg(\frac{2x}{(2k+1)\pi}\bigg)^{2}\bigg)\label{40} \\
&=\prod_{k=0}^{\infty}\bigg(1-\bigg(\frac{2x}{(2k+1)\pi}\bigg)^{4}\bigg).\nonumber	
\end{align}
From \eqref{40}, we have
\begin{align}
\tanh x-\tan x&=\frac{d}{dx}\Big(\log\cosh x+\log\cos x\Big)=\frac{d}{dx}\log\Big(\cosh x\cos x\Big) \label{41} \\
&=\frac{d}{dx}\log\bigg(\prod_{k=0}^{\infty}\bigg(1-\bigg(\frac{2x}{(2k+1)\pi}\bigg)^{4}\bigg)=\frac{d}{dx}\sum_{k=0}^{\infty}\log\bigg(1-\bigg(\frac{2x}{(2k+1)\pi}\bigg)^{4}\bigg) \nonumber \\
&=\frac{d}{dx}\sum_{k=0}^{\infty}\bigg[\log\bigg(1+\bigg(\frac{2x}{(2k+1)\pi}\bigg)^{2}\bigg)+\log\bigg(1-\bigg(\frac{2x}{(2k+1)\pi}\bigg)^{2}\bigg)\bigg] \nonumber \\ &=\sum_{k=0}^{\infty}\bigg(\cfrac{\frac{8x}{\big((2k+1)\pi\big)^{2}}}{1+\big(\frac{2x}{(2k+1)\pi}\big)^{2}}-\cfrac{\frac{8x}{\big((2k+1)\pi \big)^{2}}}{1-\big(\frac{2x}{(2k+1)\pi}\big)^{2}}\bigg)\nonumber \\
&=8x\sum_{k=0}^{\infty}\bigg(\frac{1}{(2k+1)^{2}\pi^{2}+4x^{2}}-\sum_{k=0}^{\infty}\frac{1}{(2k+1)^{2}\pi^{2}-4x^{2}}\bigg).\nonumber
\end{align}
Thus, by \eqref{41}, we get
\begin{equation}
\frac{\tanh x-\tan x}{8x}=\sum_{k=0}^{\infty}\bigg(\frac{1}{(2k+1)^{2}\pi^{2}+4x^{2}}-\frac{1}{(2k+1)^{2}\pi^{2}-4x^{2}}\bigg). \label{42}	
\end{equation}
Letting $x=\frac{\pi}{2(2n+1)},\ (n\ge 1)$ in \eqref{42}, we have
\begin{align}
&\cfrac{\tanh\big(\frac{\pi}{2(2n+1)}\big)-\tan\big(\frac{\pi}{2(2n+1)}\big)}{\frac{4\pi}{2n+1}}\nonumber\\&=\sum_{k=0}^{\infty}\bigg(\frac{1}{(2k+1)^{2}\pi^{2}+\big(\frac{\pi}{2n+1}\big)^{2}}-\frac{1}{(2k+1)^{2}\pi^{2}-\big(\frac{\pi}{2n+1}\big)^{2}}\bigg).\label{43}
\end{align}
Thus, by \eqref{43}, we get
\begin{align}
&\frac{\pi}{4(2n+1)}\bigg[\tanh\bigg(\frac{\pi}{2(2n+1)}\bigg)-\tan\bigg(\frac{\pi}{2(2n+1)}\bigg)\bigg]\label{44}\\
&=\sum_{k=0}^{\infty}\bigg(\frac{1}{(2k+1)^{2}(2n+1)^{2}+1}-\frac{1}{(2k+1)^{2}(2n+1)^{2}-1}\bigg) \nonumber\\
&=\frac{1}{2(n+1)^{2}+1}-\frac{1}{(2n+1)^{2}-1}+\frac{1}{3^{2}(2n+1)^{2}+1}-\frac{1}{3^{2}(2n+1)^{2}-1}+\cdots \nonumber
\end{align} \par
From \eqref{44}, we note that
\begin{align}
&\cfrac{4(2n+1)}{\pi\big(\tanh\big(\frac{\pi}{2(2n+1)}\big)-\tan\big(\frac{\pi}{2(2n+1)}\big)\big)}-1=\cfrac{1}{\frac{\pi\big(\tanh\big(\frac{\pi}{2(2n+1)}\big)-\tan\big(\frac{\pi}{2(2n+1)}\big)\big)}{4(2n+1)}}-1\label{45} \\
&=\cfrac{1-\frac{1}{(2n+1)^{2}+1}+\frac{1}{(2n+1)^{2}-1}-\frac{1}{3^{2}(2n+1)^{2}+1}+\frac{1}{3^{2}(2n+1)^{2}-1}-\cdots}{\frac{1}{(2n+1)^{2}+1}-\frac{1}{(2n+1)^{2}-1}+\frac{1}{3^{2}(2n+1)^{2}+1}-\frac{1}{3^{2}(2n+1)^{2}-1}+\cdots} \nonumber  \\
&=(2n+1)^{2}+\cfrac{\frac{(2n+1)^{2}+1}{(2n+1)^{2}-1}-\frac{(2n+1)^{2}+1}{3^{2}(2n+1)^{2}+1}+\frac{(2n+1)^{2}+1}{3^{2}(2n+1)^{2}-1}-\cdots}{\frac{1}{(2n+1)^{2}+1}-\frac{1}{(2n+1)^{2}-1}+\frac{1}{3^{2}(2n+1)^{2}+1}-\frac{1}{3^{2}(2n+1)^{2}-1}+\cdots}\nonumber \\
&=(2n+1)^{2}+\cfrac{((2n+1)^{2}+1)^{2}}{\cfrac{1-\frac{(2n+1)^{2}+1}{(2n+1)^{2}-1}+\frac{(2n+1)^{2}+1}{3^{2}(2n+1)^{2}+1}-\frac{(2n+1)^{2}+1}{3^{2}(2n+1)^{2}-1}+\cdots}{\frac{1}{(2n+1)^{2}-1}-\frac{1}{3^{2}(2n+1)^{2}+1}+\frac{1}{3^{2}(2n+1)^{2}-1}-\cdots}}\nonumber \\
&=(2n+1)^{2}+\cfrac{((2n+1)^{2}+1)^{2}}{-2+\cfrac{((2n+1)^{2}-1)^{2}}{\cfrac{1-\frac{(2n+1)^{2}-1}{3^{2}(2n+1)^{2}+1}+\frac{(2n+1)^{2}-1}{3^{2}(2n+1)^{2}-1}-\cdots}{\frac{1}{3^{2}(2n+1)^{2}+1}-\frac{1}{3^{2}(2n+1)^{2}-1}+\cdots}}}\nonumber\\
&=(2n+1)^{2}+\cfrac{((2n+1)^{2}+1)^{2}}{-2+\cfrac{((2n+1)^{2}-1)^{2}}{8(2n+1)^{2}+2+\cfrac{(3^{2}(2n+1)^{2}+1)^{2}}{\cfrac{1-\frac{3^{2}(2n+1)^{2}+1}{3^{2}(2n+1)^{2}-1}+\frac{3^{2}(2n+1)^{2}+1}{5^{2}(2n+1)^{2}+1}-\cdots}{\frac{1}{3^{2}(2n+1)^{2}-1}-\frac{1}{5^{2}(2n+1)^{2}+1}+\cdots}}}}\nonumber\\
&=\cdots \nonumber\\
&=(2n+1)^{2}\nonumber \\
&+\cfrac{\big((2n+1)^{2}+1\big)^{2}}{-2+\cfrac{\big((2n+1)^{2}-1\big)^{2}}{8(2n+1)^{2}+2+\cfrac{\big(3^{2}(2n+1)^{2}+1\big)^{2}}{-2+\cfrac{\big(3^{2}(2n+1)^{2}-1\big)^{2}}{16(2n+1)^{2}+2+\cfrac{(5^{2}(2n+1)^{2}+1)^{2}}{-2+\cfrac{(5^{2}(2n+1)^{2}-1)^{2}}{24(2n+1)^{2}+2+\cdots}}} }  }  }\nonumber
\end{align}
Therefore, by \eqref{45}, we obtain the following theorem.
\begin{theorem}
For $n\ge 1$, we have
\begin{align*}
&\bigg(\frac{\pi}{4}\frac{\tanh(\frac{\pi}{2(2n+1)})-\tan(\frac{\pi}{2(2n+1)})}{2n+1}\bigg)^{-1}=(2n+1)^{2}+1 \nonumber \\
&+\cfrac{\big((2n+1)^{2}+1\big)^{2}}{-2+\cfrac{\big((2n+1)^{2}-1\big)^{2}}{8(2n+1)^{2}+2+\cfrac{\big(3^{2}(2n+1)^{2}+1\big)^{2}}{-2+\cfrac{\big(3^{2}(2n+1)^{2}-1\big)^{2}}{16(2n+1)^{2}+2+\cfrac{(5^{2}(2n+1)^{2}+1)^{2}}{-2+\cfrac{(5^{2}(2n+1)^{2}-1)^{2}}{24(2n+1)^{2}+2+\cdots}}} }  }  }\nonumber	
\end{align*}
\end{theorem}
Observing that $\tanh x=1-\frac{2}{e^{2x}+1}=-\sum_{n=1}^{\infty}E_{2n-1}\frac{2^{2n-1}x^{2n-1}}{(2n-1)!}$, we have
\begin{equation}
-\tanh x=-\frac{d}{dx}\log\big(\cosh x\big)=\sum_{n=1}^{\infty}E_{2n-1}\frac{2^{2n-1}x^{2n-1}}{(2n-1)!}.\label{48}
\end{equation}
From \eqref{48} and \eqref{10}, we note that
\begin{align}
-\log\cosh x&=\sum_{n=1}^{\infty}E_{2n-1}\frac{2^{2n-1}}{(2n-1)!2n}x^{2n}
\label{49} \\
&=\sum_{n=1}^{\infty}\bigg(\frac{2x}{\pi}\bigg)^{2n}\frac{(-1)^{n}}{4n}(-1)^{n}\pi^{2n}\frac{E_{2n-1}}{(2n-1)!}\nonumber\\
&=\sum_{n=1}^{\infty}\bigg(\frac{2x}{\pi}\bigg)^{2n}\frac{(-1)^{n}}{n}\sum_{k=0}^{\infty}\frac{1}{(2k+1)^{2n}}\nonumber\\
&=\sum_{n=1}^{\infty}\bigg(\frac{2x}{\pi}\bigg)^{2n}\frac{(-1)^{n}}{n}\bigg(1-\frac{1}{4^{n}}\bigg)\zeta(2n). \nonumber
\end{align}
By \eqref{39}, we get
\begin{align}
\log\cosh(x)&=\sum_{k=0}^{\infty}\log\bigg(1+\bigg(\frac{2x}{(2k+1)\pi}\bigg)^{2}\bigg)\label{50} \\
&=\sum_{k=1}^{\infty}\log\bigg(1+\Big(\frac{2x}{(2k-1)\pi}\Big)^{2}\bigg). \nonumber
\end{align}
Therefore, by \eqref{48}, \eqref{49} and \eqref{50}, we obtain the following theorem.
\begin{theorem}
For $x\in\mathbb{R}$, we have
\begin{align*}
\log\cosh x&=\sum_{k=1}^{\infty}\log\bigg(1+\bigg(\frac{2x}{(2k-1)\pi}\bigg)^{2}\bigg)=-\sum_{k=1}^{\infty}E_{2k-1}\frac{2^{2k-1}}{(2k-1)!2k}x^{2k}\nonumber\\
&=\sum_{k=1}^{\infty}\bigg(\frac{2x}{\pi}\bigg)^{2k}\frac{(-1)^{k-1}}{k}\bigg(1-\frac{1}{4^{k}}\bigg)\zeta(2k).
\end{align*}
In particular, we have
\begin{displaymath}
\tanh x=-\sum_{k=1}^{\infty}E_{2k-1}\frac{2^{2k-1}}{(2k-1)!}x^{2k-1}.
\end{displaymath}
\end{theorem}
\section{Conclusion}
In this paper, we demonstrated that a series for $\tan x$, which gives a product expansion of the cosine function, is very useful in deriving various results. Indeed, by using this series we obtained some continued fractions, certain identities involving derivatives of $\tan x$, several expressions for $\log \cosh x$ and an identity for $\pi^{2}$. \par
Especially, we obtained two continued fractions  of $\frac{\pi}{4}$, the one in Theorem 2.1 and the other in Theorem 2.6. We remark here that the one in Theorem 2.6 gives an infinite family of continued fractions for $\frac{\pi}{4}$, since we can choose $n$ as any positive integer. \par
The authors have studied many stuffs related to Bernoulli and Euler polynomials and numbers by using many different ideas and methods. We would like to continue to carry out researches centered around Bernoulli and Euler polynomials and numbers.

\end{document}